\documentclass[11pt]{article}
\usepackage{amsmath,amsfonts,latexsym,amssymb,amscd, graphicx,pictexwd, mathrsfs, dsfont, color,mathtools}
\usepackage[title]{appendix}

\renewcommand{\P}{{\mathbb P}}
\newcommand{\E}{{\mathbb E}}

\newcommand{\bL}{{\mathbb L}}

\newcommand{\bE}{{\mathbf E}}
\newcommand{\bP}{{\mathbf P}}

\newcommand{\R}{{\mathbb R}}

\newcommand{\cF}{{\mathcal F}}

\newcommand{\cN}{{\mathcal N}}

\newcommand{\1}{{\mathds{1}}}

\newcommand{\Wass}{{\rm Wass}}
\newcommand{\Lip}{{\rm Lip}}

\newcommand{\Var}{{\bf Var\,}}
\newcommand{\Cov}{{\bf Cov}}

\newtheorem{thm}{Theorem}

\newtheorem{lem}{Lemma}

\numberwithin{equation}{section}

\providecommand{\keywords}[1]
{
	\small	
	\textbf{\textit{Keywords---}} #1
}

\begin{document}
\title{On the two-point function of the one-dimensional KPZ equation}
\author{Sergio I. L\'opez and Leandro P. R. Pimentel}
\date{\today}
\maketitle

\begin{abstract}  
In this short communication we show that basic tools from Malliavin calculus can be applied to derive the two-point function of the slope of the one-dimensional KPZ equation, starting from a two-sided Brownian motion with an arbitrary diffusion parameter, in terms of the polymer end-point annealed distribution associated to the stochastic heat equation. We also prove that this distribution is given in terms of the derivative of the variance of the solution of the KPZ equation.   
\end{abstract}

\section{Introduction and Results}
The one-dimensional Kardar-Parisi-Zhang (KPZ) equation describes the growth of an interface $h=(h_t(x)\,,\,x\in\R)_{t\geq 0}$, modeled by a height function $h_t(x)\in \R$ at position $x\in \R$ and time $t\geq 0$, whose evolution satisfies the stochastic partial differential equation 
$$\partial_t h=\frac{1}{2}\partial_x^2 h+\frac{1}{2}(\partial_x h)^2+\xi\,,$$
where $\xi$ denotes a space-time white noise. The interest in this equation started when physicists \cite{FoNeSt, KPZ} predicted that the solution to the KPZ equation has a dynamic scaling exponent equals to $3/2$, in the sense that if one rescales space by $\epsilon^{-1}$ and $h$ diffusively with $\epsilon^{1/2}$, then one has to rescale time by $\epsilon^{-3/2}$ to see non trivial fluctuations behavior as $\epsilon\downarrow 0$. Since then, the study of these fluctuations and its universality class has been a subject of extensive work \cite{Co}. Another interesting aspect of this equation is that its solutions are expected to be locally Brownian in space, and the understanding of its non-linearity brought important developments in the theory of stochastic partial differential equations \cite{Ha}. Throughout this paper we assume that $h_0(x)=\sigma W(x)$, where $(W(x)\,,\,x\in\R)$ is a standard two-sided Brownian motion, with $\sigma >0$ arbitrary, and that the white noise and the Brownian motion are independent random elements defined on some probability space $(\Omega,\cF,\bP)$, with expectation, variance and covariance denoted by $\bE$, $\Var$ and $\Cov$, respectively.

\begin{figure}\label{Fig:poly}
	\begin{center}
	\includegraphics[height=5cm]{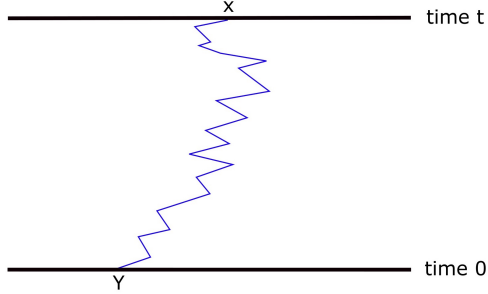}
	\caption{The time-reversal polymer realisations can be seen as paths starting at position $x$ at time $t$ and ending at the random endpoint $Y$ at time zero.}
\end{center}
\end{figure}

The KPZ equation is related to the stochastic heat equation (SHE)
\begin{equation}\label{SHE}
\partial_t Z = \frac{1}{2}\partial_x^2 Z+Z\xi \,, 
\end{equation}
by the so called Cole-Hopf transform \cite{BeGi}: $h_t(x)=\log Z_t(x)$. The Green function $Z_t(x,y)$ associated to the stochastic heat equation is the solution of \eqref{SHE} with the initial condition $Z_0(x,y)=\delta(x-y)$, the Dirac delta function for fixed $y\in\R$. The expression 
\begin{equation}\label{GreenSHE}
Z_t(x)=\int_\R Z_t(x,y)e^{\sigma W(y)}dy\,,
\end{equation}
allows us to see $Z_t(x)$ as a partition function associated to the (random) quenched density $p_{x,t}$ of a polymer endpoint $Y$, starting at position $x\in\R$ and time $t\geq 0$, and running backwards in time \cite{AlKhQu}: 
\begin{equation}\label{quenched}
p_{x,t}(y)=\frac{Z_t(x,y)e^{\sigma W(y)}}{Z_t(x)}\,\,\mbox{ and }\,\,E_{x,t} \left[f(Y)\right]=\int_\R f(y)p_{x,t}(y)dy,
\end{equation} 
see Figure \ref{Fig:poly}.

The annealed probability measure of the polymer endpoint $Y$ is defined as
\begin{equation}\label{end-point}
\P_{x,t}\left[Y\in A\right]=\bE\left[\int_A p_{x,t}(y)dy\right]\,\,\mbox{ and }\,\,\E_{x,t} \left[f(Y)\right]=\bE\left[E_{x,t} \left[f(Y)\right] \right]\,.
\end{equation}
For simplicity, when $x=0$ we remove the index from the notation. We note that $Y$ is symmetric under $\P_t$ and that 
\begin{equation}\label{stat}
\P_{x,t}\left[Y\in A\right]=\P_{t}\left[x+Y\in A\right] \,.
\end{equation}

In the study of the fluctuations of $h_t(x)=\log Z_t(x)$, the scaling function  
\begin{equation}\label{KPZscaling}
g_t(x)=\Var[h_t(x)]\,,
\end{equation}
plays a central rule, as noted in \cite{DunGuKo, GuKo}, to study KPZ equation evolving on a large torus and on the line, respectively. We recall that the dependence on the parameter $\sigma>0$ is hidden but plays an important rule in the analysis of the behavior of $h_t$. In the case $\sigma=1$ we have time stationarity, in the sense that $(h_t(x)-h_t(0)\,,\,x\in\R)\stackrel{dist.}{=}(h_0(x)\,,\,x\in\R)$, and it has been proven (Theorem $1.2$ and Proposition $1.6$ in \cite{BaQuSe}, Proposition $2.1$ and Remark $2.2$ in \cite{GuKo}) that 
\begin{equation}
g_t(x) = \mathbb{E}_{t} |Y| \sim t^{2/3}\,\, \mbox{ and } \,\, \bE\left[\partial_xh_t(z)\partial_x h_0(y)\right]= \bE\left[p_{t}(z-y)\right]=\frac{g''_t(z-y)}{2}.
\end{equation}
We also address the reader to \cite{Le,MaTh} for similar non-rigorous results in the physics literature. The left-hand side of second equality displayed above is the commonly named KPZ two-point correlation function, and this equality should be understood in the formal level, meaning that   
$$\bE\left[X_t(\phi_2)X_0(\phi_1)\right]=\int_{\R^2}\phi_2(x)\phi_1(y) \bE\left[p_{t}(x-y)\right] dxdy\,,$$
where 
$$X_t(\phi)\equiv\int_\R\phi(z)\partial_x h_t(z)dz=-\int_\R\phi'(z)h_t(z)dz\,,$$
is defined for a smooth test function $\phi:\R\to\R$ of bounded support. \\

This work can be seen as a companion communication of \cite{GuKo, Pi1}, by independently obtaining results in terms of the random polymer endpoint, including the non-stationary case, and studying the decorrelation between the initial condition and the configuration of the process as time grows to infinity. We follow the idea of bringing basic tools from Malliavin calculus (started in \cite{Pi1} in the context of the KPZ fixed point), to show directly that 
$$\bE\left[p_{t}(y)\right]=\frac{g_t''(y)}{2\sigma^2}\,\,\mbox{ and that }\,\,\bE\left[\partial_xh_t(z)\partial_x h_0(y)\right]=\frac{g_t''(z-y)}{2}\,\,\mbox{for any $\sigma>0$}\,,$$
where the second equality is understood in the formal level. The ideas contained in the proof strongly rely on the polymer representation of the solution of the KPZ equation, which also indicates that similar results can be proved mutatis-mutandis to other polymer models, as in the case of the O'Connell-Yor polymer \cite{OCYo,SeVa}, or other SHE with different multiplicative space-time noise.

By using a simple adaptation of Stein-Malliavin method for normal approximations \cite{Nu}, as developed in \cite{Pi1}, we will show that $g''_t(y)$ governs not only the decay of correlations but also the decay of the Wasserstein distance between the joint law $\theta_t$ of $(X_0^{\phi_1},X_t^{\phi_2})$ and the product measure $\eta_t$ induced by the marginals of $\theta_t$:
$$\Wass\left(\theta_t,\eta_t\right)\leq \frac{1}{\sigma\|\phi_1\|}\sqrt{\frac{\pi}{2}} \, \int_{\R^2} |\phi'_2(x)||\psi_1(y)| \frac{g_t''(x-y)}{2} dx dy\,,$$
where $\|\phi_1\|$ is the $\bL^2(\R)$ norm of $\phi_1$, $\psi'(x)=\phi_1(x)$ and $\psi(0)=0$. The Wasserstein distance between two probability measures defined on $\R^2$ is given by 
$$\Wass(\theta,\eta)=\sup_{\ell\in\Lip_1}\left\{\Big|\,\int_{\R^2} \ell\, d\theta-\int_{\R^2} \ell\, d\eta\,\Big|\right\}\,,$$
where $\Lip_1$ is the collection of all functions $\ell:\R^2\to\R$ with Lipschitz constant bounded by $1$. 

The limiting fluctuations of $h_t$ are described by the KPZ fixed point \cite{MaQuRe} for a wide class of initial profiles \cite{AmCoQu,BaQuSe,BoCoFeVe,QuSa,Vi}, which indicates that $g_t(x)\sim t^{2/3} g(x t^{-2/3})$, where $g$ is an universal scaling function given exactly by the variance of the KPZ fixed point starting from a two-sided Brownian motion (with a new diffusion coefficient depending on $\sigma$). As a consequence, one also expects that $g''_t(x)\sim t^{-2/3}g''(x t^{-2/3})$. Another interesting observation is that the quenched end-point polymer density $p_{t}(y)=p^\sigma_{t}(y)$ also makes sense for $\sigma=0$, the so called flat initial condition (see also \cite{DaZh} for localization results). From our formula for $\sigma>0$ we also get that 
$$\bE\left[p^0_{t}(y)\right]=\lim_{\sigma\searrow 0}\bE\left[p^\sigma_{t}(y)\right]=\lim_{\sigma\searrow 0}\frac{g_t''(y)}{2\sigma^2}\,.$$   
(Recall that $g_t(x)=g^\sigma_t(x)$ also depends on $\sigma>0$.) 

Our main results are summarized in the theorem below, in a slightly different but equivalent formulation that will fit the proofs better. We denote 
$$\phi_2\star\phi_1(u)=\int_{ \mathbb{R} } \phi_2(z)\phi_1(z+u)dz\,,$$ 
the cross-correlation between the test functions $\phi_2$ and $\phi_1$. Note that (by the symmetry of $g''_t$)
$$\int_{\R^2}\phi_2(x)\phi_1(y)\frac{g_t''(x-y)}{2}dxdy=\int_\R\phi_2\star\phi_1(u)\frac{g_t''(u)}{2}du\,.$$
   
\begin{thm}\label{Main}
Recall the annealed end-point probability measure \eqref{end-point} and the scaling function \eqref{KPZscaling}. We have that 
\begin{enumerate}
\item $g'_t(y)=\sigma^2\left(2\P_{t}\left[Y\leq y\right]-1\right)$;
\item $\bE\left[X_t(\phi_2)X_0(\phi_1)\right]=\sigma^2\E_{t} \left[\phi_2\star\phi_1(Y)\right]$;
\item $\Wass\left(\theta_t,\eta_t\right)\leq \frac{\sigma}{\|\phi_1\|}\sqrt{\frac{\pi}{2}}  \E_{t}\left[|\phi'_2|\star|\psi_1|\left(Y\right)\right]$.
\end{enumerate}
\end{thm} 

\section{Proofs}

The Malliavin derivative, with respect to the standard two-sided Brownian motion $W$, of a random variable $X$ defined on $(\Omega,\cF,\bP)$ is a $\bL^2(\Omega\times \R)$ valued random element denoted by $(DX(u)\,,\,u\in\R)$ \cite{Nu}. For instance, for a smooth test function $\phi:\R\to\R$, we have that 
\begin{equation}\label{MalDer}
D\left(\int_\R\phi(z) dW(z)\right)(u)=\phi(u)\,.
\end{equation}
One of the main tools  in Malliavin calculus is the integration by parts formula, which states that 
\begin{equation}\label{IP}
\bE\left[X\int_\R\phi(z) dW(z) \right]=\bE\left[\langle DX,\phi\rangle\right]\,,
\end{equation} 
where $\langle\cdot,\cdot\rangle$ is the usual inner product in $\bL^2(\R)$.   
Let us define the function 
\begin{equation}\label{zeta}
u\in\R\,\mapsto\,\1_y(u)=\left\{\begin{array}{ll}\,\,\,\,\1_{(0,y]}(u)& \mbox{ if } y>0\,,\\
\,\,\,\,0& \mbox{ if } y=0\,,\\
- \1_{(y,0]}(u) &\mbox{ if } y<0\,.\,\end{array}\right.
\end{equation}
that will be used several times in the following, as well as the usual notation $x\vee y=\max\{x,y\}$, $x\wedge y=\min\{x,y\}$, $x^+=x\vee 0$ and $x^-=-(x\wedge 0)$.

The Malliavin derivative \eqref{MalDer} of $\sigma W(y)=\int_\R \sigma \1_y(u) dW(u)$ is given by 
\begin{equation*}
D\left(\sigma W(y)\right)(u)=\sigma\1_y(u)\,,
\end{equation*}
and by the chain rule, one gets that
$$D(e^{\sigma W(y)})(u)=e^{\sigma W(y)}D\left(\sigma W(y)\right)(u)= \sigma e^{\sigma W(y)}\1_y(u)\,.$$
Using \eqref{GreenSHE}, one has that $Z_t(x)$ is Malliavin differentiable with respect to $W$ \cite{GuKo}. By switching the derivative with the Lebesgue integral in \eqref{GreenSHE} (Proposition $3.4.3$ in \cite{Nu2}), and using the previous calculation, one obtains that (notice that $Z_t(x,y)$ does not depend on $W$)
$$D\left(Z_t(x)\right)(u)=\sigma\int_{\R}Z_t(x,y)e^{\sigma W(y)}\1_y(u)dy\,,$$
and by the chain rule again,
$$
D\left(h_t(x)\right)(u)=D\left(\log Z_t(x)\right)(u)=\frac{1}{Z_t(x)}D\left(Z_t(x)\right)(u)=\sigma\frac{\int_{\R}Z_t(x,y)e^{\sigma W(y)}\1_y(u)dy}{Z_t(x)}\,. 
$$
Writing this expression in terms of the polymer end-point density \eqref{quenched}, one finally concludes that
\begin{equation}\label{DerKPZ}
D\left(h_t(x)\right)(u)=\sigma E_{x,t}\left[\1_{Y}(u)\right]\,. 
\end{equation}

\begin{lem}\label{Variance}
We have that 
$$g_t(x)=\Var\left[h_t(x)\right]=\Var\left[h_t(0)\right]-\sigma^2 |x|+2\sigma^2 G_t(x)\,,$$ 
where,
$$G_t(x):=\left\{\begin{array}{ll}\E_{t}\left[\left(x\wedge(x+Y)\right)^+\right] & \mbox{ if } x\geq 0\,,\\
 \E_{t}\left[\left(x\vee(x+Y)\right)^-\right]&\mbox{ if } x\leq 0\,.\end{array}\right.$$ 
\end{lem}

\noindent\paragraph{\bf Proof Lemma \ref{Variance}}
Use \eqref{GreenSHE}, together with shift invariance of the Brownian motion and translation invariance of the white noise, to see that 
$$h_t(x)-h_0(x)=\log\int_\R Z_t(x,y)e^{\sigma(W(y)-W(x))}dy\stackrel{dist.}{=}\log\int_\R Z_t(0,y)e^{\sigma W(y)}dy=h_t(0)\,.$$
Since, $\E\left[h_0(x)\right]=\sigma\E\left[W(x)\right]=0$, we get that 
\begin{eqnarray*}
\Var\left[h_t(x)\right]&=&\Var\left[h_t(x)-h_0(x)\right]+\Var\left[h_0(x)\right]+2\Cov\left[h_t(x)-h_0(x),h_0(x)\right]\\
&=&\Var\left[h_t(0)\right]-\Var\left[h_0(x)\right]+2\Cov\left[h_t(x),h_0(x)\right]\\
&=&\Var\left[h_t(0)\right]-\sigma^2 |x|+2\bE\left[h_t(x)h_0(x)\right]\,.
\end{eqnarray*}
Integration by parts \eqref{IP} from Malliavin calculus, combined with \eqref{DerKPZ}, implies that 
$$\bE\left[h_t(x)h_0(x)\right] =\sigma\bE\left[h_t(x)\int_\R\1_x(z)dW(z)\right]=\sigma\bE\left[\langle Dh_t(x)\,,\,\1_x\rangle\right]=\sigma^2\bE\left[\langle E_{x,t}\left[\1_{Y}(\cdot)\right]\,,\,\1_x\rangle\right]\,.$$
By Fubini's theorem, for $x\geq 0$,
$$\langle E_{x,t}\left[\1_{Y}(\cdot)\right]\,,\,\1_x\rangle=\int_0^xE_{x,t}\left[\1_{Y}(u)\right]du=E_{x,t}\left[\int_0^x\1_{Y}(u)du\right]=E_{x,t}\left[\left(x\wedge Y\right)^+\right]\,,$$
while for $x\leq 0$,
$$ \langle E_{x,t}\left[\1_{Y}(\cdot)\right]\,,\,\1_x\rangle=-\int_x^0E_{x,t}\left[\1_{Y}(u)\right]du=E_{x,t}\left[-\int_x^0\1_{Y}(u)du\right]=E_{x,t}\left[\left(x\vee Y\right)^-\right]\,.$$
To conclude the proof, use \eqref{end-point} and \eqref{stat}.
\hfill$\Box$\\

\begin{lem}\label{CrossFor}
$$\bE\left[\langle DX_t(\phi_2),\phi_1\rangle\right]=\sigma\E_{t}\left[\phi_2\star\phi_1(Y)\right]\,\mbox{ and }\,\bE\left[\left|\langle DX_t(\phi_2),\phi_1\rangle\right|\right]\leq\sigma\E_{t}\left[|\phi'_2|\star|\psi_1|(Y)\right]\,.$$
\end{lem}

\noindent\paragraph{\bf Proof Lemma \ref{CrossFor}} 
Recall \eqref{zeta} and 
$$\psi_1(x)=\int_\R\phi_1(u)\1_x(u)du\,\mbox{ (the primitive of $\phi_1$)}\,.$$ 
By interchanging Malliavin derivative and integral and using \eqref{DerKPZ},
$$DX_t^{\sigma}(\phi_2)(u)=-\sigma\int_\R \phi'_2(z)E_{z,t}\left[\1_{Y}(u)\right]dz\,,$$
and hence, by Fubini's theorem,
\begin{eqnarray}
\nonumber\langle DX_t(\phi_2),\phi_1\rangle&=&\int_\R DX_t(\phi_2)(u)\phi_1(u)du\\
\nonumber&=&-\sigma\int_\R\left(\int_\R\phi'_2(z)E_{z,t}\left[\1_{Y} (u)\right]\phi_1(u)dz\right)du\\
\label{Fub}&=&-\sigma\int_\R\phi'_2(z)E_{z,t}\left[\psi_1 (Y)\right]dz\,.
\end{eqnarray}
On one hand, Fubini's theorem again, \eqref{Fub}, \eqref{end-point} and \eqref{stat},
$$\bE\left[\langle DX_t(\phi_2),\phi_1\rangle\right]=-\sigma\int_\R\phi'_2(z)\E_{t}\left[\psi_1(z+Y)\right]dz\,.$$
On the other hand, by standard integration by parts,
$$\int_\R\phi'_2(z)\psi_1(z+Y)dz=-\int_\R\phi_2(z)\phi_1(z+Y)dz\,,$$
and hence, by Fubini's theorem once again, 
\begin{eqnarray*}
\int_\R\phi'_2(z)\E_{t}\left[\psi_1(z+Y)\right]dz&=& \E_{t}\left[\int_\R\phi'_2(z)\psi_1(z+Y)dz \right]\\
&=& -\E_{t}\left[\int_\R\phi_2(z)\phi_1(z+Y)dz \right]\,,
\end{eqnarray*}
and we can conclude that
$$ \bE\left[\langle DX_t(\phi_2),\phi_1\rangle\right]=\sigma\E_{t}\left[\phi_2\star\phi_1(Y)\right]\,.$$
Taking the absolute value into the integrals, a similar calculation leads to 
$$\bE\left[\left|\langle DX_t(\phi_2),\phi_1\rangle\right|\right]\leq\sigma\E_{t}\left[|\phi'_2|\star|\psi_1|(Y)\right]\,.$$
\hfill$\Box$\\

\noindent\paragraph{\bf Proof Theorem \ref{Main}} 
To prove that 
$$g'_t(y)=\sigma^2\left(2\P_{t}\left[Y\leq y\right]-1\right)\,,$$
one uses the formula obtained in Lemma \ref{Variance}, and notice that, for $x\geq 0$,
\begin{eqnarray*}
\E_{t}\left[\left(x\wedge(x+Y)\right)^+\right]&=& \E_{t}\left[(x+Y)\1_{\{Y\in(-x,0]\}}+x\1_{\{Y>0]\}}\right]\\
&=&\E_{t}\left[Y\1_{\{Y\in(-x,0]\}}\right]+x\P_{t}\left[Y>-x\right]\,.
\end{eqnarray*}
and that, for $x\leq 0$, 
\begin{eqnarray*}
\E_{t}\left[\left(x\vee(x+Y)\right)^-\right]&=& \E_{t}\left[-(x+Y)\1_{\{Y\in(0,-x]\}}-x\1_{\{Y\leq 0]\}}\right]\\
&=&-\E_{t}\left[Y\1_{\{Y\in(0,-x]\}}\right]-x\P_{t}\left[Y\leq -x\right]\,.
\end{eqnarray*}
Thus, for $x\geq 0$, 
$$G'_t(x)=\P_{t}\left[Y>-x\right]=\P_{t}\left[Y\leq x\right]\,,$$
and hence 
$$g'_t(x)= -\sigma^2+2\sigma^2\P_{t}\left[Y\leq x\right]=\sigma^2\left(2\P_{t}\left[Y\leq x\right]-1\right)\,,$$
while for $x\leq 0$, 
$$ G'_t(x)=-\P_{t}\left[Y\leq -x\right]=\P_{t}\left[Y\leq x\right]-1\,,$$
and hence 
$$g'_t(x)= \sigma^2+2\sigma^2\left(\P_{t}\left[Y\leq x\right]-1\right)=\sigma^2\left(2\P_{t}\left[Y\leq x\right]-1\right)\,.$$

By integration by parts \eqref{IP},
$$\bE\left[X_t(\phi_2)X_0(\phi_1)\right]=\sigma\bE\left[X_t(\phi_2)\int_\R\phi_1(z)dW(z)\right]=\sigma\bE\left[\langle DX_t(\phi_2),\phi_1\rangle\right]\,,$$
and the second item of Theorem \ref{Main} follows from Lemma \ref{CrossFor}.

To prove the upper bound on the Wasserstein distance between $\theta_t$ and $\eta_t$ we give a brief explanation of the Stein-Malliavin method for asymptotic independence, as developed in \cite{Pi1}.  Let $\sigma_1>0$ and $X_1\sim N(0,\sigma_1^2)$ be a normal random variable with mean zero and variance $\sigma_1^2$. Define the functional operator $\cN$, acting on differentiable functions $f:\R^2\to\R$, by
\begin{equation}\label{SteinChar}
\cN f(x_1,x_2):=\sigma_1^2\partial_{x_1} f(x_1,x_2)- x_1 f(x_1,x_2)\,.
\end{equation}
Given a random variable $X_2$ that is independent of $X_1$, the operator $\cN$ characterizes the joint product measure $\eta$ of $(X_1,X_2)$ in the sense that $X_1\sim N(0,\sigma_1^2)$ and $X_2$ is independent of $X_1$ if and only if $\bE\left[\cN f(X_1,X_2)\right]=0$ for all nice functions $f$. In this way, given a random vector $(X_1,X_2)$ with probability measure $\theta$, and such that $X_1\sim N(0,\sigma_1^2)$, we might estimate the distance to independence by looking at $\bE\left[\cN f(X_1,X_2)\right]$. This idea can be made precise by using the unique bounded solution $f_\ell$ of the partial differential equation
$$ \cN f(x_1,x_2)= \ell(x_1,x_2)-\bE\left[\ell(X_1,x_2)\right]\,,$$ 
in such way that 
$$\bE\left[\cN f_\ell(X_1,X_2)\right]=\int_{\R^2} \ell\, d\theta-\int_{\R^2} \ell\, d\eta\,\implies\,\Wass(\theta,\eta)\leq \sup_{\ell}\left| \bE\left[\cN f_\ell(X_1,X_2)\right]\right|\,.$$
In our context, we have $X_1=X_0(\phi_1)=\sigma \int_\R\phi_1(z)dW(z)\sim N(0,\sigma_1^2)$ with $\sigma_1^2=\sigma^2\|\phi_1\|^2$, and $X_2=X_t(\phi_2)$. By using that
$$Df(X_1,X_2)=\partial_{x_1}f(X_1,X_2)DX_1+\partial_{x_2} f(X_1,X_2) DX_2\,,$$
and integration by parts \eqref{IP}, we have that 
\begin{multline*}
\bE\left[X_1f(X_1,X_2)\right]=\sigma\bE\left[\langle Df(X_1,X_2),\phi_1\rangle\right]\\
=\sigma^2\|\phi_1\|^2\bE\left[\partial_{x_1}f(X_1,X_2)\right]+\sigma\bE\left[ \partial_{x_2} f(X_1,X_2) \langle DX_2,\phi_1\rangle\right] \,,
\end{multline*}
and hence
$$\bE\left[\cN f(X_1,X_2)\right]= -\sigma\bE\left[ \partial_{x_2} f(X_1,X_2) \langle DX_2,\phi_1\rangle\right]\,.$$
By Lemma 8 \cite{Pi1}, 
$$\sup_{\R^2}\left|\partial_{x_2} f_\ell(x_1,x_2)\right|\leq\frac{1}{\sigma_1}\sqrt{\frac{\pi}{2}}\sup_{\R^2}\left|\partial_{x_2} \ell(x_1,x_2)\right|\leq\frac{1}{\sigma_1}\sqrt{\frac{\pi}{2}}\,$$
(recall that $\ell$ has Lipschitz constant bounded by $1$), which finally yields the upper bound
\begin{multline*}
\Wass(\theta_t,\eta_t)\leq \sup_{\ell}\left| \bE\left[\cN f_\ell(X_1,X_2)\right]\right|\\
\leq \frac{1}{\|\phi_1\|}\sqrt{\frac{\pi}{2}}\bE\left[\left|\langle DX_2,\phi_1\rangle\right|\right]\leq \frac{\sigma}{\|\phi_1\|}\sqrt{\frac{\pi}{2}}\E_{0,t}\left[ |\phi'_2|\star|\psi_1|(Y)\right]\,,
\end{multline*}
where in the last equality we use Lemma \ref{CrossFor}.

\hfill$\Box$\\

\textbf{Acknowledgements}
Sergio I. L\'opez was supported by DGAPA-UNAM grants PAPIIT IN-102822 and PASPA 2021-II. Leandro P. R. Pimentel was supported by the National Council of Scientific Researches (CNPQ, Brazil) Grant 305356/2019-4.


\begin{thebibliography}{10}

\bibitem{AlKhQu}
\textsc{Alberts, T., Khanin, K. and Quastel, J.} (2014). The Continuum Directed Random Polymer. \textit{J. Stat. Phys.} {\bf  154}, 305--326. 

\bibitem{AmCoQu} 
\textsc{Amir, G., Corwin, I., and Quastel, J.} (2011). Probability distribution of the free energy of the continuum directed random polymer in 1 + 1 dimensions. \textit{Comm. Pure Appl. Math.} {\bf 64},  466--537.

\bibitem{BaQuSe}
\textsc{Balázs, M., Quastel, J., and Seppäläinen, T.} (2011). Fluctuation exponent of the KPZ/stochastic Burgers equation. \textit{J. Amer. Math. Soc.} {\bf 24}, 683--708.

\bibitem{BeGi}
\textsc{Bertini, L., and Giacomin, G.} (1997). Stochastic Burgers and KPZ equations from particle systems. \textit{Comm. Math. Phys.} {\bf 183}, 571--607.

\bibitem{BoCoFeVe}
\textsc{Borodin, A., Corwin, I., Ferrari, P. L., and V\"eto, B.} (2015). Height fluctuations for the stationary KPZ equation. \textit{Math. Phys. Anal. Geom.} {\bf 18}, 20. 

\bibitem{Co}
\textsc{Corwin, I.} (2012). The Kardar-Parisi-Zhang equation and universality class. \textit{Random Matrices: Theory Appl.} {\bf 1}, 1130001.

\bibitem{DaZh}
\textsc{Das, S. and Zhu, W.} (2022). Localization of the continuum directed random polymer. Available at https://arxiv.org/abs/2203.03607. 

\bibitem{DunGuKo}
\textsc{Dunlap, A. and Gu, Y., and Komorowski, T.} (2021). Fluctuations of the KPZ equation on a large torus. To appear in Comm. Pure and Appl. Math., available at https://arxiv.org/abs/2111.03650.


\bibitem{FoNeSt}
\textsc{Forster, D., Nelson, D.R., and Stephen, M.J.} (1977). Large-distance and long-time properties of a randomly stirred fluid. \textit{Phys. Rev. A} {\bf 16}, 732--749.

\bibitem{GuKo}
\textsc{Gu, Y., and Komorowski, T.} (2023). Another look at Bal\'azs-Quastel-Sepp\"al\"ainen theorem. \textit{Trans. Amer. Math. Soc.}  {\bf 376}, 2947--2962.

\bibitem{Ha}
\textsc{Hairer, M.} (2013). Solving the KPZ equation. \textit{Ann. Math.} {\bf 178}, 559--664.

\bibitem{KPZ}
\textsc{Kardar, M., Parisi, G., and Zhang, Y.Z.} (1986). Dynamic scaling of growing interfaces. \textit{Phys. Rev. Lett.} {\bf 56}, 889--892.

\bibitem{Le}
\textsc{Le Doussal, P.} (2017). Maximum of an Airy process plus Brownian motion and memory in Kardar-Parisi-Zhang growth. \textit{Phys. Rev. E} \textbf{96} 060101.

\bibitem{MaTh}
\textsc{Maes, C.} and \textsc{Thiery, T.} (2017). Midpoint Distribution of Directed Polymers in the Stationary Regime: Exact Result Through Linear Response. \textit{J. Stat. Phys.} \textbf{168} 937--963.

\bibitem{MaQuRe}
\textsc{K. Matetski, K., Quastel, J., and Remenik, D.} (2021). The KPZ fixed point.
\textit{Acta Math.} {\bf 227}, 115--203.

\bibitem{Nu}
\textsc{Nualart, D.} (2019). Malliavin calculus and normal approximations. \textit{Ensaios Matem\'aticos} {\bf 34}, 1--74.

\bibitem{Nu2}
\textsc{Nualart, D., Nualart, E.} (2018). Introduction to Malliavin calculus. \textit{IMS textbooks}, Cambridge University Press.

\bibitem{OCYo}
\textsc{O’Connell, N., Yor, M.} (2001). Brownian Analogues of Burke’s Theorem. \textit{Stoch. Process 
 Appl.} {\bf 96} 285--304.

\bibitem{Pi1}
\textsc{Pimentel, L.P.R.} (2020). Integration by parts and the KPZ two-point function. \textit{Ann. Probab.} {\bf 50}(5), 1755--1780.

\bibitem{QuSa} 
\textsc{Quastel, J., and Sarkar, S.} (2023). Convergence of exclusion processes and the KPZ equation to the KPZ fixed point. \textit{J. Amer. Math. Soc.} {\bf 36} 251--289. 

\bibitem{SeVa}
\textsc{Sepp\"al\"ainen, T., Valk\'o, B.} (2010). Bounds for scaling exponents for a 1+1 dimensional
directed polymer in a Brownian environment. \textit{ALEA} {\bf 7}, 451--476.

\bibitem{Vi}
\textsc{Virag, B.} (2020). The heat and the landscape I. Available at https://arxiv.org/pdf/2008.07241.


\end{thebibliography}
\end{document}